\documentclass[12pt,twoside]{article}
\usepackage{amsmath,amsthm,amssymb}
\usepackage{times}
\usepackage{enumerate}
\usepackage{delarray}
\usepackage{xcolor,graphicx}
\usepackage{pgfplots}
\usepackage{tikz,pgf}

\theoremstyle{definition}
\newtheorem{theorem}{Theorem}[section]
\newtheorem{remark}[theorem]{Remark}

\newtheorem{proposition}[theorem]{Proposition}

\numberwithin{equation}{section}
\newcommand{\subjclass}[1]{\bigskip\noindent\emph{2010 Mathematics Subject Classification:}\enspace#1}
\newcommand{\keywords}[1]{\noindent\emph{Keywords:}\enspace#1}

\newcommand{\om}{\omega}

\renewcommand{\th}{\theta}

\newcommand{\cU}{\mathcal{U}}
\newcommand{\cS}{\mathcal{S}_\Phi}
\newcommand{\cZ}{\mathcal{Z}}

\newcommand\enne{\mathbf{N}}

\newcommand\erre{\mathbf{R}}

\newcommand{\pa}{\partial}

\newcommand{\ov}{\overline}
\newcommand{\dyle}{\displaystyle}
\newcommand{\dimo}{{\bf Proof.}\quad}
\newcommand{\finedim}{{\unskip\nobreak\hfil\penalty50
  \hskip2em\hbox{}\nobreak\hfil\mbox{\rule{1ex}{1ex} \qquad}
  \parfillskip=0pt\finalhyphendemerits=0\par\medskip}}

\begin{document}

\baselineskip=17pt


\title{Some remarks on segregation of \(k\) species\\ 
in strongly competing systems}
\author{Flavia Lanzara \\
Dipartimento di Matematica {\it G. Castelnuovo}, \\
{\it Sapienza} Universit\`a di Roma, \\
p.le A. Moro 5 00185 Roma \\
flavia.lanzara@uniroma1.it \\
Eugenio Montefusco \\
Dipartimento di Matematica {\it G. Castelnuovo}, \\ 
{\it Sapienza} Universit\`a di Roma, \\
p.le A. Moro 5 00185 Roma \\
eugenio.montefusco@uniroma1.it}

\date{}

\maketitle

\begin{abstract}
Spatial segregation occurs in population dynamics when \(k\) 
species interact in a highly competitive way. As a model for the 
study of this phenomenon, we consider the competition-diffusion 
system of \(k\) differential equations
\[
-\Delta u_i(x)=-\mu u_i (x)\dyle\sum_{j\neq i} u_j (x) \quad i=1,...,k
\]
in a domain \(D\) with appropriate boundary conditions. Any \(u_i\) 
represents a population density and the parameter \(\mu\) determines 
the interaction strength between the populations. The purpose of this 
paper is to study the geometry of the limiting configuration as 
\(\mu\longrightarrow+\infty\) on a planar domain for any number of 
species. If \(k\) is even we show that some limiting configurations 
are strictly connected to the solution of a Dirichlet problem for the 
Laplace equation. 

\subjclass{Primary 35Bxx, 35J47; Secondary 92D25.}

\keywords{Pattern formation, Spatial segregation, Strong competition}
\end{abstract}

\section{Introduction and setting of the problem}

When two or more species live in proximity and share the same basic 
requirements, they usually compete for resources, habitat or territory. 
Only the strongest prevails, driving the weaker competitors to 
extinction. This is the principle of competitive exclusion (also 
known as Gause's law). One species wins because its members are more 
efficient at finding resources, which leads to an increase 
in population. This means that a population of competitors finds less 
of the same resources and cannot grow at its maximal capacity \cite{lek}.

According to the competitive exclusion principle, many 
competing species cannot coexist under very strong competition,
but when spatial movements are permitted more than one species can 
coexist thanks to the segregation of their habitats. For a theoretical
discussion and some experimental results see \cite{gau,ha}.

From a mathematical viewpoint the determination of the configuration 
of the habitat segregation for some populations is an interesting 
problem which can be modelled by an optimal (in a suitable sense) 
partition of a domain; for example in the papers 
\cite{ctv1,ctv2,ctv3,ctv4,tt12} the problem is studied modelling the 
interspecies competition with a large interaction term in an elliptic 
system of partial differential equations inspired by classical models
in populations dynamics. In \cite{bz,bz2,dwz} the problem is modelled 
as a Cauchy problem for a parabolic system of semilinear partial 
differential equations describing the dynamics of the densities 
of different species. 
In the evolutive case, in particular see \cite{dwz}, it is proved that 
some populations can vanish under the competition of other species;
moreover, in \cite{bz,bz2} the authors are able to estimate the number 
of the long-term surviving populations and other interesting qualitative 
properties of the spatial distributions of interacting populations.
Note that also the study of the territoriality, that is how different 
groups of the same species divide an area, avoiding to effectively 
fight for resources, can be viewed as a habitat segregation produced
by competition (see, for example, \cite{bdf08,halw,mk05}); moreover, 
this kind of competition is a struggle between competitors having 
the same features, that is between perfect competitors.

Others approaches have been considered. In \cite{cdnp}, for example, 
a phase segregation problem is studied by minimization of integral 
functionals. The authors obtain results showing mixing or separation
of the two phases, depending of the strength of the interaction between 
the two species (that is depending on some relations between the
parameters in the model). The main difference with our study is that, in  
\cite{cdnp}, the densities have no to satisfy a boundary datum, so that 
the optimal partition of the domain is related only to the geometry of 
the domain and the interaction between the parameters.

As a model for the study of the segregation  phenomena, we consider the 
competition-diffusion system of \(k\) differential equations
\begin{equation}\label{diffeqk}
\begin{array}\{{cl}.
-\Delta u_i(x)=-\mu u_i (x)\dyle\sum_{j\neq i} u_j (x)& \hbox{ in }D\,, \\
u_i(x)\geq 0 & \hbox{ in }D\,, \\
u_i(x)=\phi_i(x) & \hbox{ on }\pa D\,.
\end{array} \qquad i=1,...,k\,.
\end{equation}
Here  \(D\subseteq\erre^n\) is  an open bounded, simply connected  
domain with smooth boundary $\pa D$. In this paper we consider the 
case \(n=2\). We assume that the function \(\Phi=(\phi_1,...,\phi_k)\) 
is an {\it admissible datum} that is  \(\phi_i\in W^{1,\infty}(\pa D)\),   
\(\phi_i\geq 0\), \(i=1,...,k\), \(\phi_i\cdot\phi_j=0\) a.e. in 
\(\partial D\) for \(i\neq j\), the sets \(\{\phi_i>0\}\) are nonempty, 
open connected arcs and the function \(\sum_{i=1}^k\phi_i\) vanishes 
exactly in \(k\) points of \(\partial D\) (the endpoints of the 
\(\phi_i\)'s supports).

The system  \eqref{diffeqk} governs the steady states of \(k\) competing 
species coexisting in the same domain \(D\). Any \(u_i\) represents a 
population density and the parameter {\(\mu>0\)} determines the interaction 
strength between the populations. In this model the competition between 
two species is depicted without direct reference to the resources they 
share, rather, it is assumed that the presence of each population leads 
to a depression of its competitor's  growth rate.

If \(\Phi\) is admissible, the existence of positive solutions of 
\eqref{diffeqk} for any positive \(\mu\) is proved in \cite{ctv2} using 
Leray-Schauder degree theory. The uniqueness is proved in \cite{wz}, 
using the sub- and super-solution method.

Let define the class of segregated densities
\[
\cU_\Phi = \begin{array}\{{rl}\}
U=(u_1,...,u_k)\in (H^1(D))^k \hbox{ : }  & u_i=\phi_i \hbox{ on }\pa D \\
& u_i\geq 0 \hbox{ in } D \\
 & u_i\cdot u_j=0 \hbox{ for \(i\neq j\) a.e. in } D \\
\end{array}
\]
and the class
\[
\cS=\begin{array}\{{rl}\}
U=(u_1,...,u_k)\in \cU_\Phi\hbox{ : } 
 & -\Delta u_i\leq 0 \hbox{ in } D \\
 & -\Delta \big(u_i-\sum_{j\neq i} u_j \big)\geq 0 \hbox{ in } D \\
\end{array}\,.
\]
Let \(U^{(\mu)}=(u_{1,\mu},...,u_{k,\mu})\) be the solution of 
\eqref{diffeqk} for every \(\mu>0\). In \cite{ctv2} it is proved 
that there exists \(\overline{U}=\{\bar u_1,...,\bar u_k\}\in\cU\) 
such that, up to subsequences, \(u_{i,\mu}\longrightarrow\bar{u}_i\) 
in \(H^1(D)\) and \(\cS\) contains all the asymptotic limits of 
\eqref{diffeqk} that is \(\overline{U}\in \cS\). 

The uniqueness of the limit solution of \eqref{diffeqk} as 
\(\mu\to+\infty\) was proved in \cite{ctv2} in the case \(k=2\) 
and in \cite{ctv4} in the case of \(k=3\) and in dimension \(n=2\). 
Specifically, the authors prove that the class \(\cS\) consists 
of one element. In \cite{wz} it is proved that \(\cS\) consists of 
one element also in the case of arbitrary dimension and arbitrary 
number of species. A different proof of uniqueness of the limit 
configuration, based on the maximum principle and on the 
qualitative properties of the elements of \(\cS\), is given in \cite{ab}.  

The description of the qualitative properties of the limiting 
configurations in the planar case (i.e. \(n=2\))  was considered 
in \cite{ctv4} for \(k=3\) and in \cite{lm19} for \(k=4\). 
The aim of this paper is to describe the geometry of the limiting 
configurations in the planar case, for any number of species. 

The outline of the paper is as follows. In Section 2 we recall some 
basic facts and known results  of the class \(\cS\), which will be 
used later. In Section 3 we study  the geometry of the limiting 
configuration in any number of species. If \(k=2s\), starting from 
the argument used in \cite{lm19} for \(k=4\) species, we 
prove that some limiting configurations are strictly related to the 
solution of a Dirichlet problem for the Laplace equation. 
Our results rely on the construction of a harmonic function which  
assumes the value \(\sum_{j=1}^{2s}(-1)^j\phi_j\) on \(\partial D\). 
This function has an even number of nodal regions compatible with 
an alternate sign rule. We emphasize that this construction cannot 
be done in the case of odd \(k\).
This will be the object of a forthcoming paper.
In Section 4 we focus on the case of \(k=6\) number of species.

\section{Basic facts}

In this section we recall some basic facts that will play an 
important role in our study. Suppose that \(D\) is a simply 
connected domain in \(\erre^2\). Due to the conformal 
invariance of the problem, with no loss of generality 
we can assume
\[
D=B(O,1)=\{x=(x_1,x_2) \in\erre^2: |x|<1\}
\]
and consider the class
\[
\cS=\begin{array}\{{rl}\}
&U=(u_1,...,u_k)\in (H^1(D))^k \hbox{ : } 
u_i\geq 0 \hbox{ in } D ,\>  u_i=\phi_i \hbox{ on }\pa D\\
&u_i\cdot u_j=0 \hbox{ for } i\neq j,  -\Delta u_i\leq 0,   
-\Delta \big(u_i-\sum_{j\neq i} u_j \big)\geq 0 \hbox{ in } D 
\end{array}.
\]
The study of \(\cS\) provides the understanding of the 
segregated states of \(k\) species induced by strong 
competition. If  $\Phi=(\phi_1,...,\phi_k)$ is an admissible 
datum then \(\phi_i, i=1,...,k,\) are positive in their 
supports, the sets  $\{\phi_i>0\}\subset \partial D$ are 
open connected arcs  and $\Phi=\sum_{i=1}^k \phi_i$ vanishes 
at exactly $k$ points of $\pa D$, the endpoints 
$p_1,....,p_k$ in counter clockwise order. 

In the following we will denote by \(U\) both the \(k-\)uple 
\((u_1,...,u_k)\) and the function \(\sum_{i=1}^k u_i\). 
For any $U\in\cS$ define the {\it nodal regions}
\[
\om_i=\{p\in D:u_i(p)>0\} \qquad i=1,...,k,
\]
the {\it multiplicity of a point} $p \in\ov{D}$ with 
respect to $U$:
\[
m(p)=\#\{i: |\om_i\cap B_r(p)|>0 \quad \forall r>0\}
\]
where $B_r(p)=\{q\in\erre^2: |p-q|<r\}$, and the {\it interfaces} 
between two densities
\[
\Gamma_{ij}=\pa \om_i\cap \pa \om_j \cap \{p\in D: m(p)=2\}.
\]
The element \(U\in\cS\) defines exactly  \(k\) nodal regions and
\[
\overline D=\overline{\om}_1\cup...\cup \overline{\om}_k.
\]

Let us summarize the basic 
properties of the elements $U\in\cS$:
\begin{itemize}
\item[(s1)] each $u_i\in W^{1,\infty}(\ov{D})$ (\cite[Theorem 8.4]{ctv3}). 
It follows that $u_i\in C(\ov{D})$, $\om_i$ is open and \(p\in\om_i\) 
implies $m(p)=1$. By standard regularity theory for elliptic 
equations we also have that \(u_i\in C^\infty(\omega_i)\);
\item[(s2)] each $\om_i$ is connected and each $\Gamma_{ij}$ 
is either empty or a connected arc starting from a point \(p_i\in\partial D\)  
(\cite[Remark 2.1]{ctv4});
\item[(s3)] $u_i$ is harmonic in $\om_i$; $u_i-u_j$ is harmonic 
in $D\setminus \cup_{h\neq i,j}\overline \om_h$ 
\cite[Proposition 6.3]{ctv3};
\item[(s4)] if $p\in D$ satisfies \(m(p)=2\), then 
(\cite[Remark 6.4 ]{ctv3})
\[
\lim_{\om_i\ni q\to p} \nabla u_i(q)
=-\lim_{\om_j\ni q\to p} \nabla u_j(q) ;
\]
\item[(s5)] $U\in W^{1,\infty}(\ov{D})$ and if $p\in D$ satisfies 
\(m(p)=2\), then 
\[
|\nabla U(p)| = \lim_{\om_i\ni q\to p} |\nabla u_i(q)|
=\lim_{\om_j\ni q\to p} |\nabla u_j(q)| \neq 0
\] 
and the set  $\{q:m(q)=2\}$ is locally a $C^1$-curve through 
$p$ ending either at points with higher multiplicity, 
or at the boundary $\pa D$  \cite[Lemma 9.4]{ctv3};
\item[(s6)] if $p\in D$ satisfies $m(p)\geq 3$, then 
$|\nabla U(q)|\to 0$, as $q\to p$ \cite[Theorem 9.3]{ctv3};
\item[(s7)] the set $\{p\in D: m(p)\geq 3\}$ consists of a finite 
number of points \cite[Lemma 9.11]{ctv3};
\item[(s8)] if $p\in D$ with $m(p)=h\geq 3$ then there 
exists \(\th_0\in (-\pi,\pi]\) such that
\begin{equation}\label{local}
U(r,\th)=r^{h/2} \left| \cos\left( \frac{h}{2} 
(\th +\th_0)\right) \right|+o(r^{h/2})
\end{equation}
as $r\to 0$, where $(r,\th)$ is a system of polar coordinates 
around \(p\) \cite[Theorem  9.6]{ctv3}.
\end{itemize}

\begin{remark} The asymptotic formula \eqref{local} describes 
the behavior of \(U\) in a neighborhood of a multiple point in \(D\). 
As a consequence, at multiple point \(U\in\cS\) shares the 
angle in equal parts. This property does not hold true if 
\(U\in\cS\) has a multiple point \(p\) on the boundary 
\(\partial D\). 
\end{remark}

\section{Results on \(k\) species}

Let \(U\in\cS\), we define the set of points of multiplicity 
greater than or equal to $h\in\enne$
\[
\cZ_h(U)=\{p\in \overline D: m(p)\geq h\}.
\]
The set $\cZ_h(U)$ consists of a finite number of isolated 
points \cite[Lemma 9.11 and Theorem 9.13]{ctv3}.

\begin{proposition}\label{Pk.1}
Let \(k\geq 3\) and \(U\in\cS\), then \(\cZ_3(U)\) is nonempty 
and does not contain points of multiplicity higher than \(k\).
\end{proposition}

\dimo If \(\cZ_3(U)=\emptyset\) then the interfaces 
\(\Gamma_{ij}\) between any two densities do not intersect in 
\(\overline D\). The function \(\Phi=\sum_{i=1}^k\phi_i\) vanishes 
in  exactly \(k\) points and any \(\Gamma_{ij}\) links a  point 
\(p_i\in \partial D\) to a point \(p_j\in \partial D\), \(i\neq j\). 
Therefore, if \(k\) is odd  there exists at least a point 
\(p_\ell\in\partial D\) which belongs to two interfaces. 
Then \(p_\ell\in\cZ_3(U)\) and we get a contradiction.

If \(k\) is even, since \(\phi\) vanishes in  
exactly \(k\) points, then there are only \(k/2+1\) interfaces 
which are nonempty, with endpoints on the boundary and do not 
intersect. This contradicts the fact that \(U\in\cS\) defines 
\(k\) nodal regions. It also implies that \(\cZ_s(U)=\emptyset\) 
for \(s>k\).
\finedim

Let \(U\in\cS\). We associate to any \(p\in\cZ_3(U)\) the number
\[
i(p)=m(p)-2.
\]

\begin{proposition}\label{pr:formula}
Let \(k\geq 2\). The following relation holds
\begin{equation}\label{formula}
k-2=\sum_{p\in\cZ_3(U)}i(p).
\end{equation}
\end{proposition}

\dimo If \(k=2\) then \(\cZ_3(U)=\emptyset\) and \eqref{formula} 
is trivially satisfied. Let \(k\geq 3\). First of all we want 
to point out that the set 
\(\left(\cup\overline{\Gamma}_{ij}\right)\cup \partial D\), that 
is the union of the interfaces between any  two species and the 
boundary of the disk forms a planar, connected graph, whose vertices 
are the points in \(\cZ_3(U)\) and the zeros of the boundary datum. 
From classical arguments in graph theory (essentially the Euler 
polyhedral formula, see for example \cite[Theorem 1.5.2]{j}) 
it follows that 
\[
n-m+f=2
\]
where \(n\) is the number of the vertices of the graph, \(m\) the 
number of the edges and \(f\) the number of the faces. In our case 
we have that \(n =\sharp\{\cZ_3(U)\cap D\}+k\), that is the multiple 
points in the disk plus the \(k\) zeros of  the boundary datum,
\(m=\sharp\{\Gamma_{ij}\}+k\), the number of the arcs 
\(\Gamma_{ij}\neq\emptyset\) plus \(k\) (the number of the arcs 
\(\{\phi_i>0\}\subseteq \partial D\)), and \(f\) is the number 
of the nodal regions \(\omega_i\) plus \(1\)
so that  \(f=k+1\). Then, it follows that 
\[
m=n+f-2=2k+\sharp\{\cZ_3(U)\cap D\}-1.
\]
Now we want to point out that, for any \(p\in\cZ_3(U)\cap D\), 
the number \(m(p)\) corresponds to the number of the arcs 
\(\Gamma_{ij}\) such that \(p\in\overline{\Gamma}_{ij}\), whileas 
for the others vertices \(p\in\partial D\) (the zeros of \(\Phi\) 
on the boundary) it holds that \(m(p)\) is the number of the arcs 
\(\Gamma_{ij}\) such that \(p\in\overline{\Gamma}_{ij}\) augmented 
of \(1\). Then, recalling that for any vertex 	
\(p\in\partial D\setminus\cZ_3(U)\) we have \(m(p)=2\), we can write that
\[
\begin{aligned}
\sum_{p\in\cZ_3(U)} i(p)
	& = \sum_{p\in\cZ_3(U)} \left[m(p)-2\right] 
		= \sum_{p\in\cZ_3(U)\cap D} \left[m(p)-2\right] 
		+ \sum_{p \textrm{ vertices on }\partial D} \left[m(p)-2\right] \\
	& = \sum_{p\in\cZ_3(U)\cap D} 
		\left[\sharp\{\Gamma_{ij}: p\in \ov{\Gamma}_{ij}\}-2\right] 
		+ \sum_{p \textrm{ vertices on }\partial D} 
		\left[\sharp\{\Gamma_{ij}: p\in \ov{\Gamma}_{ij}\}+1-2\right] \\
	& = \sum_{\substack{p \textrm{ vertices}\\ \textrm{of the graph}}}
		\left[\sharp\{\Gamma_{ij}: p\in \ov{\Gamma}_{ij}\}-2\right] 
		+ k
\end{aligned}
\]
since the vertices on the boundary are exactly \(k\). Then, it follows that
\[
\begin{aligned}
\sum_{p\in\cZ_3(U)} i(p)
	& = \sum_{\substack{p \textrm{ vertices}\\ \textrm{of the graph}}} 
		\left[\sharp\{\Gamma_{ij}: p\in \ov{\Gamma}_{ij}\}-2\right] + k 
		= \sum_{\substack{p \textrm{ vertices}\\ \textrm{of the graph}}} 
		\left[\sharp\{\Gamma_{ij}: p\in \ov{\Gamma}_{ij}\}\right] -2n + k \\
	& = 2(m-k) -2n +k 
		= 2(f-2) -k = k-2
\end{aligned}
\]
since the edges are the union of the interfaces \(\Gamma_{ij}\) 
and of the \(k\) arcs \(\{\phi_i>0\}\subseteq \partial D\), and 
summing over the vertices any edge is counted twice.
\finedim

\begin{remark}
Note that, following \cite{j}, the set 
\(\Gamma=\cup{\ov{\Gamma}}_{ij}\), that is the union of the 
interfaces between two species,  the zeros of the boundary 
datum \(\Phi\) and the points in \(\cZ_3(U)\) form a tree, since 
it is a planar, acyclic, connected graph. The zeros of \(\Phi\) 
are the leaves and the multiple points in \(\cZ_3(U)\) are the 
other vertices. From classical arguments in graph theory it 
follows that the number of the arcs composing \(\Gamma\) is 
\((k+\#\{\cZ_3(U)\}-1)\).
\end{remark}

\begin{remark}
Let \(k\geq 3\). We want to point out that identity 
\eqref{formula} implies that  
\[
1\leq \#\{\cZ_3(U)\}\leq k-2\,.
\]
\end{remark}

In the following we assume that the number of species   
\(k\) is even, that is \(k=2\,s\), \(s>1\). 
Then we can define a harmonic function having opposite 
signs on adjacent nodal regions and strictly connected 
to \(U\in\cS\).

Consider the boundary value problem
\begin{equation}\label{armonico}
\begin{array}\{{cl}.
-\Delta\psi= 0 & \hbox{ in }D \\
	\psi=\phi & \hbox{ on }\pa D \,.
\end{array}
\end{equation}
\begin{proposition}
Let \(\psi_a\) be the solution of  \eqref{armonico} with 
boundary datum \(\Phi^a= \sum_{j=1}^{2s} (-1)^j\phi_j\) 
and let \(U\in\cS\). \\
1. If \(U=|\psi_a|\) then for any \(p\in \cZ_3(U)\cap D\) 
we have that \(m(p)\) is even.\\
2. If there exists \(p\in \cZ_3(U)\cap D\) such that \(m(p)\) 
is odd then \(U\neq |\psi_a|\). 
\end{proposition}

\dimo
1. Let  \(U=|\psi_a|\) and \(p\in \cZ_3(U)\). 
If \(p\in D\) then \(U(p)=0\) and \(\nabla U(p)=(0,0)\)  
(\cite[Theorem 9.3]{ctv3}). It follows that \(p\) is a 
critical point for \(\psi_a\) at level \(0\). By standard 
theory of harmonic functions the zero set of \(\psi_a\) around 
a critical point at level 0 is made by (at least) \(4\) 
half-lines, meeting with equal angles. 
We infer that locally around \(p\) the function \(\psi_a\) 
defines \(q\) nodal components with \(q\geq 4\) and \(q\) 
is even because \(\psi_a\) has alternate positive or 
negative sign on adjacent sets. 

2. Let \(p\in \cZ_3(U)\) such that \(m(p)=q\) with odd 
\(q\). If \(p\in D\) and \(U = |\psi_a|\) then \(p\) is 
a critical point for \(\psi_a\)  at level \(0\). 
We infer that locally around \(p\) the function 
\(\psi_a\) defines \(q\) nodal region, and \(q\) is odd. 
This is a contradiction because \(\psi_a\) has 
alternate positive or negative sign on adjacent sets. 
\finedim

\begin{theorem}\label{aleharmonic}
Let \(\Phi=(\phi_1,...,\phi_{2s})\) be an admissible datum.
The harmonic function \(\psi_a\) which solves \eqref{armonico}
with boundary datum \(\Phi^a= \sum_{j=1}^{2s} (-1)^j\phi_j\) 
possesses at most \(s-1\) critical points \(q\) in \(D\) 
such that \(\psi_a(q)=0\).
\end{theorem}
 \dimo Since \(\Phi\) is an admissible datum, the solution of 
\eqref{armonico} with boundary datum \(\Phi^a\) vanishes at 
exactly \(2s\) points on \(\partial D\), each arc 
\(\{\phi_j>0\}\subset\partial D\) is connected 
and \(\psi_a\) has different signs on adjacent arcs.

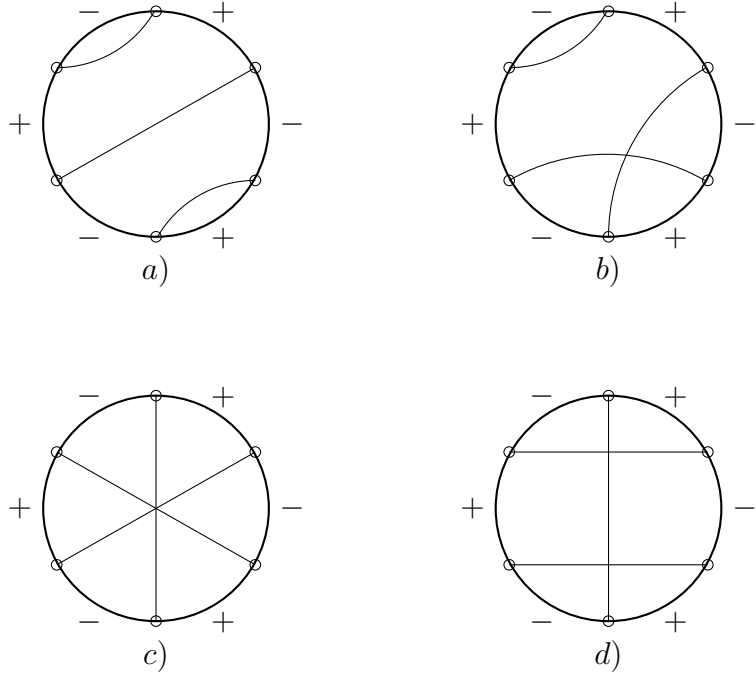
\begin{figure}[h!]
\begin{center}
\begin{tikzpicture}
\draw[thick] (0,0) circle (1.5cm);
\draw (-1.32,-0.75) -- (1.32,0.75);
\draw (-1.32,0.75) arc (270:330:1.5);
\draw (1.32,-0.75) arc (90:150:1.5);
\draw (0,-1.5) circle (2pt);
\draw (0,1.5) circle (2pt);
\draw (-1.32,-0.75) circle (2pt);
\draw (1.32,0.75) circle (2pt);
\draw (1.32,-0.75) circle (2pt);
\draw (-1.32,0.75) circle (2pt);
\node [anchor= south east] at (1.2, 1.2) {\(+\)};
\node [anchor= south west] at (-1.2, 1.2) {\(-\)};
\node [anchor= east] at (-1.5, 0) {\(+\)};
\node [anchor= west] at (1.5, 0) {\(-\)};
\node [anchor= south west] at (-1.2, -1.8) {\(-\)};
\node [anchor= south east] at (1.2, -1.8) {\(+\)};
\node [anchor= north] at (0, -1.6) {\(a)\)};
\end{tikzpicture}
\qquad \qquad
\begin{tikzpicture}
\draw[thick] (0,0) circle (1.5cm);
\draw (0,-1.5) arc (180:120:2.6);
\draw (-1.32,-0.75) arc (120:60:2.6);
\draw (-1.32,0.75) arc (270:330:1.5);
\draw (0,-1.5) circle (2pt);
\draw (0,1.5) circle (2pt);
\draw (-1.32,-0.75) circle (2pt);
\draw (1.32,0.75) circle (2pt);
\draw (1.32,-0.75) circle (2pt);
\draw (-1.32,0.75) circle (2pt);
\node [anchor= south east] at (1.2, 1.2) {\(+\)};
\node [anchor= south west] at (-1.2, 1.2) {\(-\)};
\node [anchor= east] at (-1.5, 0) {\(+\)};
\node [anchor= west] at (1.5, 0) {\(-\)};
\node [anchor= south west] at (-1.2, -1.8) {\(-\)};
\node [anchor= south east] at (1.2, -1.8) {\(+\)};
\node [anchor= north] at (0, -1.6) {\(b)\)};
\end{tikzpicture}
\\
\vspace{1cm}
\begin{tikzpicture}
\draw[thick] (0,0) circle (1.5cm);
\draw (0,-1.5) -- (0,1.5);
\draw (-1.32,-0.75) -- (1.32,0.75);
\draw (1.32,-0.75) -- (-1.32,0.75);
\draw (0,-1.5) circle (2pt);
\draw (0,1.5) circle (2pt);
\draw (-1.32,-0.75) circle (2pt);
\draw (1.32,0.75) circle (2pt);
\draw (1.32,-0.75) circle (2pt);
\draw (-1.32,0.75) circle (2pt);
\node [anchor= south east] at (1.2, 1.2) {\(+\)};
\node [anchor= south west] at (-1.2, 1.2) {\(-\)};
\node [anchor= east] at (-1.5, 0) {\(+\)};
\node [anchor= west] at (1.5, 0) {\(-\)};
\node [anchor= south west] at (-1.2, -1.8) {\(-\)};
\node [anchor= south east] at (1.2, -1.8) {\(+\)};
\node [anchor= north] at (0, -1.6) {\(c)\)};
\end{tikzpicture}
\qquad \qquad 
\begin{tikzpicture}
\draw[thick] (0,0) circle (1.5cm);
\draw (0,1.5) -- (0,-1.5);
\draw (-1.32,-0.75) -- (1.32,-0.75);
\draw (1.32,0.75) -- (-1.32,0.75);
\draw (1.32,0.75) circle (2pt);
\draw (1.32,-0.75) circle (2pt);
\draw (-1.32,0.75) circle (2pt);
\draw (-1.32,-0.75) circle (2pt);
\draw (0,1.5) circle (2pt);
\draw (0,-1.5) circle (2pt);
\node [anchor= south east] at (1.2, 1.2) {\(+\)};
\node [anchor= south west] at (-1.2, 1.2) {\(-\)};
\node [anchor= east] at (-1.5, 0) {\(+\)};
\node [anchor= west] at (1.5, 0) {\(-\)};
\node [anchor= south west] at (-1.2, -1.8) {\(-\)};
\node [anchor= south east] at (1.2, -1.8) {\(+\)};
\node [anchor= north] at (0, -1.6) {\(d)\)};
\end{tikzpicture}
\end{center}
\caption{The level  set \(\Gamma=\{x\in D:\psi_a(x)=0\}\) and 
the nodal components of \(\psi_a\), where the function \(\psi_a\) 
solves \eqref{armonico} with boundary datum 
\(\Phi^a= \sum_{j=1}^6 (-1)^j\phi_j\). Three different situations 
can occur: the function \(\psi_a\) has \(4\) nodal regions 
(figure \(a)\)), \(5\) nodal regions (figure \(b)\)) or \(6\) 
nodal regions (figures \(c)\) and \(d)\)). }\label{figtypebis}
\end{figure}

Since a harmonic function does not admit closed level lines, 
the set \(\Gamma=\{x\in D:\psi_a(x)=0\}\) has no closed loop. 
We infer that \(\psi_a\) has alternate positive or negative 
sign on \(\rho\) sets, with \(s+1\leq \rho\leq 2s\): the nodal 
components of \(\psi_a\) (see figure \ref{figtypebis} 
in the case \(k=6\)).
The zero set of the harmonic function \(\psi_a\) around a 
critical point at level 0 is made by (at least) \(4\) half-lines, 
meeting with equal angles. We infer that locally around each 
critical point at level \(0\) the function \(\psi_a\) defines 
\(\ell\) nodal components with \(4\leq \ell\leq \rho\) and \(\ell\) 
is even because \(\psi_a\) has alternate positive or negative 
sign on adjacent sets. 

Suppose that \(\psi_a\) has \(r\) critical points 
\(q_1,...,q_r\) in \(D\) such that \(\psi_a(q_j)=0\), 
\(j=1,...,r\), \(r\geq 1\). Then \(q_j\in\cZ_3(|\psi_a|)\) 
with \(m(q_j)\geq 4\), \(j=1,...,r\).
If \(\psi_a\) defines \(2s\) nodal regions i.e. \(\rho=2s\) 
then the function \(U=|\psi_a|\) is an 
element of \(\cS\). From \eqref{formula},
\[
2s-2=\sum_{p\in\cZ_3(U)}i(p)
\geq  \sum_{j=1}^{r}  i(q_j)=\sum_{j=1}^{r} m(q_j)-2r\geq 4r-2r=2r
\]
that is \(r\leq s-1\). If \(\psi_a\) defines \(\rho\) nodal 
regions with \(\rho<2s\)  then, by repeating the same argument 
in the proof of formula \eqref{formula}, we obtain that
\begin{equation}\label{formula2}
2\rho-2-2s=\sum_{p\in\cZ_3(|\psi_a|)}i(p).
\end{equation}
Then,
\[
2\rho-2-2s=\sum_{p\in\cZ_3(|\psi_a|)}i(p)
\geq  \sum_{j=1}^{r}  i(q_j)=\sum_{j=1}^{r} m(q_j)-2r\geq 4r-2r=2r.
\]
We infer that \(r<\rho-s-1<s-1\).
\finedim

If the solution \(\psi_a\) to \eqref{armonico} with 
alternate sign on adjacent arcs has exactly \(s-1\) 
critical points at level zero we can say much more.
\bigskip

\begin{proposition}\label{2spoint2}
Let  \(\Phi=(\phi_1,...,\phi_{2s})\) be an admissible boundary 
datum and suppose that the harmonic function \(\psi_a\), solution 
to \eqref{armonico} with boundary datum 
\(\Phi^a=\sum_{j=1}^{2s}(-1)^j \phi_j\), has 
\(q_1,...,q_{s-1}\)  critical points in \(D\) such that 
\(\psi_a(q_i)=0\), \(i=1,...,s-1\).  
Then \(q_1,...,q_{s-1}\) are \(4\)-points for the function 
\(U=|\psi_a|\in\cS\).
\end{proposition}

\dimo The zero set of the harmonic function \(\psi_a\) around a 
critical point at level 0 is made by (at least) \(4\) half-lines, 
meeting with equal angles. We infer that locally around each \(q_i\), 
\(i=1,...,s-1\),  the function \(\psi_a\) defines \(\ell\) nodal 
components with \(\ell\geq 4\). The function \(\psi_a\) defines 
\(\rho\leq 2s\) nodal regions. From formula  \eqref{formula2}, 
\[
2\rho-2-2s=\sum_{p\in\cZ_3(|\psi_a|)}(m(p)-2)
\geq \sum_{j=1}^{s-1}(m(q_j)-2)\geq 2(s-1).
\]
We infer that \(\rho\geq 2s\), hence \(\rho=2s\).

The function \(U=|\psi_a|\) is nonnegative, satisfies the 
boundary datum and has exactly  \(2s\) nodal regions. 
This function generates an element of \(\cS\), with 
boundary datum \(\Phi\) and \(\{q_1,...,q_{s-1}\}\in \cZ_3(U)\).
Hence, we have
\[
\{q_1,...,q_{s-1}\}\subseteq \cZ_3(U),\quad m(q_j)=m_j
\]
with \(m_j\geq 4\), \(j=1,...,s-1\). From \eqref{formula} we get
\[
\sum_{j=1}^{s-1}i(q_j)\leq 2s-2
\]
or equivalently
\[
\sum_{j=1}^{s-1}m_j-2(s-1)\leq 2s-2 
\iff \sum_{j=1}^{s-1}m_j\leq 4(s-1)\,.
\]
Since \(m_j\geq 4\) we deduce that \(m_j=4\), \(j=1,...,s-1\). \finedim

\begin{proposition}\label{iffharmonic} 
Let \(\Phi=(\phi_1,...,\phi_{2s})\) be an admissible datum. 
If \(U\in \cS\) possesses a  
\(2s\)-point \(a_U\) in \(\overline{D}\) then \(U=|\psi_a|\), 
where \(\psi_a\) is the harmonic function such that 
\(\psi_a=\Phi^a=\sum_{j=1}^{2s} (-1)^j\phi_j\) on 
\(\partial D\). If \(a_U\in D\) then \(a_U\) is a 
critical point for \(\psi_a\) at zero  level.
\end{proposition}

\dimo Let \(U=(u_1,...,u_{2s})\) be an element of  
\(\cS\) with a \(2s\)-point \(a_U\), then 
there exist \(2s-1\) or \(2s\) arcs connecting \(a_U\) 
to any of the isolated zeros of the boundary datum 
(note that there are \(2s-1\) arcs if and only if 
\(a_U\in\partial D\)). Then the function 
\(\psi_a=\sum_{j=1}^{2s}(-1)^ju_j\) is harmonic in 
\(D\setminus\{a_U\}\) (see \cite[Proposition 6.3]{ctv3}), 
moreover \(\psi_a\) is bounded so, by Schwarz's removable 
singularity principle (see \cite[Proposition 11.1]{p}) 
\(\psi_a\) is harmonic in \(D\) and, by construction, 
\(U=|\psi_a|\) in \(D\) and \(U=|\Phi^a|=\Phi\) on 
\(\partial D\).

If \(a_U\in\cZ_3(U)\cap D\) then, for the property (s6),  
\(U(a_U)=0\) and \(\nabla U(a_U)=(0,0)\). Hence  
\(\psi_a(a_U)=0\) and \(\nabla \psi_a(a_U)=(0,0)\).
\finedim

We can generalized Proposition \ref{iffharmonic} as follows.

\begin{proposition}\label{iffharmonic2} 
Let \(\Phi=(\phi_1,...,\phi_{2s})\) be an admissible datum. 
Let \(U\in \cS\) be such that \(\cZ_3(U)\) contains 
only points with even multiplicity in \(\overline{D}\). 
Then  \(U=|\psi_a|\), where \(\psi_a\) is the solution 
of the Dirichlet problem \eqref{armonico} with boundary 
datum \(\psi_a=\Phi^a= \sum_{j=1}^{2s} (-1)^j\phi_j\) on 
\(\partial D\). Moreover  any \(p\in\cZ_3(U)\cap D\) 
is a critical point for \(\psi_a\) at zero level.
\end{proposition}

\dimo Let \(U\) be an element of  \(\cS\) such that 
\(\cZ_3(U)=\{q_1,...,q_r\}\) with \(m(q_i)=2\ell_i\geq 4\), 
\(i=1,...,r\). Then, for every \(i=1,...,r\),  there exist
\(s_1^{(i)},...,s_{2\ell_i}^{(i)}\in \{1,...,2s\}\) such that 

i) \(\overline\om_{s^{(i)}_1}\cap...
\cap\overline\om_{s^{(i)}_{2\ell_i}}=\{q_i\}\);

ii) \(\sum_{j=1}^{2\ell_i}(-1)^ju_{s_j}\) is harmonic in 
\(\dot{\Omega}_i\setminus \{q_i\}\) where
\(\Omega_i=\overline\om_{s^{(i)}_1}\cup...
\cup\overline\om_{s^{(i)}_{2\ell_i}}\) and \(\dot{\Omega}_i\) 
denotes the interior of \(\Omega_i\).

Keeping in mind that \(u_i=0\) in \(\om_j\), \(i\neq j\), 
then \(\psi_a=\sum_{j=1}^{2s} (-1)^j u_j\) is harmonic 
in \(\dot\Omega_i\setminus \{q_i\}\), \(i=1,...,r\). 
It is easy to see that 
\(\dot\Omega_1\cup...\cup\dot\Omega_r=D\setminus\{q_1,...,q_r\}\).
Then the function \(\psi_a=\sum_{j=1}^{2s}(-1)^ju_j\) is harmonic 
in \(D\setminus\{q_1,...,q_r\}\) (see \cite[Proposition 6.3]{ctv3}), 
moreover \(\psi_a\) is bounded so, by Schwarz's removable 
singularity principle (see \cite[Proposition 11.1]{p}) 
\(\psi_a\) is harmonic in \(D\) and, by construction, 
\(U=|\psi_a|\) in \(D\) and \(U=|\Phi^a|=\Phi\) on \(\partial D\).

If \(q\in\cZ_3(U)\cap D\) then \(U(q)=0\) and \(\nabla U(q)=(0,0)\). 
Hence  \(\psi_a(q)=0\) and \(\nabla \psi_a(q)=(0,0)\).
\finedim

The next Proposition gives  conditions on the 
admissible datum \(\Phi\) such that the limiting 
configuration has a point with multiplicity \(2s\) in \(D\).

\begin{proposition}\label{CS}
Let  \(\Phi=(\phi_1,...,\phi_{2s})\) be an admissible boundary datum and let 
\(\psi_a\) be the solution to \eqref{armonico} 
with boundary datum \(\Phi^a=\sum_{j=1}^{2s}(-1)^j \phi_j\).  
Let \(p\in{D}\) such that  \(\Phi^a\) satisfies the conditions
\begin{eqnarray}\label{cond1}
\int_{\partial D} \Phi^a\left(\frac{\zeta+p}{\overline{p}\zeta+1}
\right) T_j(\zeta_1)ds_\zeta=0,\quad j=0,1,...,s-1\\
\label{cond2}
\int_{\partial D} \Phi^a\left(\frac{\zeta+p}{\overline{p}\zeta+1}
\right) \zeta_2 U_{j-1}(\zeta_1)ds_\zeta=0
,\quad j=1,...,s-1
\end{eqnarray}
with  \(\zeta=(\zeta_1,\zeta_2)\), where \(T_j\) and \(U_j\) 
denote the   Chebyshev polynomials of the first and second 
kind, respectively. Then \(p\) is a \(2s\)-point 
for the function \(U=|\psi_a|\in\cS\).
\end{proposition}

\dimo We introduce the transformation 
\begin{equation}\label{conforme}
x=R_p(\zeta)=\frac{\zeta+p}{\overline{p}\zeta+1}\,.
\end{equation}
Here we identificate the complex numbers \(x=x_1+i x_2\) 
and \(\zeta=\zeta_1+i\zeta_2\) with the points \((x_1,x_2)\) 
and \((\zeta_1,\zeta_2)\in\erre^2\), respectively. \(R_p\) 
is a conformal map which maps the unit disk \(\overline{D}\) 
into itself such that \(R_p(\partial D)=\partial D\) 
and \(R_p(0)=p\). Set
\begin{equation}\label{TR}
\widetilde \Psi_a(\zeta)=\psi_a(R_p(\zeta)),\qquad
\widetilde \Phi^a(\zeta)=\Phi^a(R_p(\zeta))\,.
\end{equation}
Then \(\widetilde \Psi_a\) solves the problem
\begin{equation}\label{DP2}
\begin{array}\{{cl}.
-\Delta\widetilde \Psi_a= 0 & \hbox{ in }D\,, \\
	\widetilde \Psi_a=\widetilde \Phi^a & \hbox{ on }\pa D\,. \\
\end{array}
\end{equation}
Introducing a system of polar coordinats \((r,\theta)\), 
we can write the Fourier expansion of \(\widetilde \Psi_a\)
\begin{equation}\label{fourier}
\widetilde \Psi_a(r,\theta)=\frac{A_0}{2}
+\sum_{j=1}^\infty (A_j \cos(j\theta)+B_j \sin(j\theta))r^j\,,
\qquad \zeta=(r,\theta)\,.
\end{equation}
Let \(T_j\) and \(U_j\) denote the Chebychev polynomials 
of the first and second kind, respectively. Keeping in mind 
the representation \(T_j(\cos(\theta))=\cos(j\theta)\) and 
\(U_{j-1}(\cos(\theta))=\sin(j\theta)/\sin(\theta)\) 
(\cite[22.3.15-16] {as}) from the conditions 
\eqref{cond1}-\eqref{cond2} we get
\begin{eqnarray}\label{cond1bis}
A_j:=\frac{1}{\pi}\int_{-\pi}^{\pi}\widetilde\Phi^a(e^{i\theta})
\cos(j\theta)d\theta=0,\quad j=0,1,...,s-1
\\\label{cond2bis}
B_j:=\frac{1}{\pi}\int_{-\pi}^{\pi}\widetilde\Phi^a(e^{i\theta})
\sin(j\theta)d\theta=0,\quad j=1,...,s-1.
\end{eqnarray}
It follows that, around the origin,
\[
\widetilde \Psi_a(r,\theta)=\sum_{j=s}^\infty (A_j \cos(j\theta)
+B_j \sin(j\theta))r^j\,,
\qquad \zeta=(r,\theta)\,.
\]
We have \((A_s,B_s)\neq (0,0)\). Indeed, if not, let 
\((A_\nu,B_\nu)\neq (0,0)\) where \(\nu>s\) is the index of the 
first nonzero Fourier component.  Then there would be \(2\nu\) 
arcs starting form the origin, on which \(\widetilde \Psi_a\)  
vanishes. Since a harmonic function does not admit closed level 
lines, this contradicts the fact that \(\widetilde \Phi^a\) 
has exactly \(2s\) zeroes.

Therefore 
\(U(x)=|\psi_a(x)|=|\widetilde \Psi_a (R^{-1}_p(x))|\)
is nonnegative, satisfies the boundary datum \(\Phi\) and 
has exactly \(2s\) nodal regions. This function generates 
an element of \(\cS\), with datum \(\Phi\) and the 
\(2s\)-point \(p\) (see also \cite[Lemma 3.2]{ctv4}).
\finedim

Conversely,

\begin{proposition}\label{CN} 
Suppose that  \(\Phi=(\phi_1,...,\phi_{2s})\) is an admissible 
datum and the function  \(U\in \cS \) generates a configuration 
with a  \(2s\)-point in \(p\in {D}\). 
Then \(\Phi^a=\sum_{j=1}^{2s}(-1)^j \phi_j\) satisfies the 
conditions \eqref{cond1}-\eqref{cond2}.
\end{proposition}

\dimo For Proposition \ref{iffharmonic}, we have 
\(U=|\psi_a|\) with  \(\psi_a\)  solution of 
\eqref{armonico} with boundary datum 
\(\Phi^a=\sum_{j=1}^{2s} (-1)^j \phi_j\). Then the 
function \(\widetilde \Psi_a(\zeta)=\psi_a(R_p(\zeta))\), 
with \(R_p\) in \eqref{conforme}, solves 
\eqref{DP2}, it is given by 
\begin{equation}\label{poisson}
\widetilde \Psi_a(\zeta)=\frac{1-|\zeta|^2}{2\pi}
\int_{\partial D} \frac{\widetilde \Phi^a (\eta)}{|\zeta-\eta|^2}
ds_\eta
\end{equation}
and it belongs to \(C^2(D)\cap C^0(\overline{D})\) 
(cf. \cite[(2.27)]{gt}).

On the other hand in the Fourier expansion \eqref{fourier} 
of \(\widetilde \Psi_a\) we have
\[
A_0=A_j=B_j=0,\qquad j=1,...,s-1
\]
and \((A_s,B_s)\neq (0,0)\).  The result follows from 
\eqref{cond1bis}-\eqref{cond2bis}.\finedim

\begin{proposition}\label{CNeS1}
Let \(\Phi=(\phi_1,...,\phi_{2s})\) be an admissible boundary 
datum, \(\Phi^a=\sum_{j=1}^{2s}(-1)^j \phi_j\) and \(p\in{D}\). 
Conditions \eqref{cond1}-\eqref{cond2} are equivalent to  
\begin{eqnarray}\label{cond1ter}
\int_{\partial D} \Phi^a\left(\frac{\zeta+p}{\overline{p}\zeta+1}
\right) \zeta_1^{j-h} \zeta_2^hds_\zeta=0,
\qquad h=0,...,j; j=0,...,s-1
\end{eqnarray}
with  \(\zeta=(\zeta_1,\zeta_2)\).
\end{proposition}

\dimo From the representations
\begin{eqnarray*}\label{sinj}
\sin(j\theta)=\sum_{\substack{{h=1}\\{h\>\> odd}}}^j (-1)^{\frac{h-1}{2}}  
\begin{pmatrix} j\\h 
\end{pmatrix} \cos^{j-h}\theta \sin^h \theta,\qquad 
\cos(j\theta)=\sum_{{\substack{{h=0}\\{h\>\> even}}}}^j (-1)^{\frac{h}{2}}  
\begin{pmatrix} j\\h 
\end{pmatrix} \cos^{j-h}\theta \sin^h \theta
\end{eqnarray*}
we infer that, setting \(\zeta_1=\cos(\theta), \zeta_2=\sin(\theta)\),
\begin{eqnarray}\label{sinj2}
\zeta_2 U_{j-1}(\zeta_1)=\sum_{\substack{{h=1}\\{h\>\> odd}}}^j  
(-1)^{\frac{h-1}{2}}  
\begin{pmatrix} 
j\\h 
\end{pmatrix} 
\zeta_1^{j-h} \zeta_2^{h} ,\quad \zeta_1^2+\zeta_2^2=1, \quad  j\geq 1 
\\ \label{cosj2}
T_j(\zeta_1)=\sum_{{\substack{{h=0}\\{h\>\> even}}}}^j
(-1)^{\frac{h}{2}}  \begin{pmatrix} j\\h 
\end{pmatrix} \zeta_1^{j-h} \zeta_2^h,\quad  \zeta_1^2+\zeta_2^2=1, j\geq 0.
\end{eqnarray}
If conditions \eqref{cond1ter} are satisfied, from the 
representations \eqref{sinj2}-\eqref{cosj2}, we deduce 
that even conditions \eqref{cond1}-\eqref{cond2} are.

Conversely, let conditions \eqref{cond1}-\eqref{cond2} hold.  
For \(j=0\), since \(T_0\equiv 1\), we get
\begin{equation}\label{j0}
\int_{\partial D} \Phi^a\left(\frac{\zeta+p}{\overline{p}\zeta+1}
\right)ds_\zeta=0.
\end{equation}
For \(j=1\), since \(T_1(\zeta_1)=\zeta_1\) and \(U_0\equiv 1\), we get
\[
\int_{\partial D} \Phi^a\left(\frac{\zeta+p}{\overline{p}\zeta+1}
\right) \zeta_1ds_\zeta=\int_{\partial D} \Phi^a
\left(\frac{\zeta+p}{\overline{p}\zeta+1}
\right) \zeta_2ds_\zeta=0.
\]
For \(j=2\), from the relations \(T_2(\zeta_1)=2 \zeta_1^2 -1\) and 
\(U_1(\zeta_1)=2 \zeta_1\), and \eqref{j0} we get
\[0=\int_{\partial D} \Phi^a\left(\frac{\zeta+p}{\overline{p}\zeta+1}
\right) (2\zeta_1^2-1) ds_\zeta
=2\int_{\partial D} \Phi^a\left(\frac{\zeta+p}{\overline{p}\zeta+1}
\right) \zeta_1^2ds_\zeta
\]
\[0=\int_{\partial D} \Phi^a\left(\frac{\zeta+p}{\overline{p}\zeta+1}
\right) \zeta_2 (2\zeta_1) ds_\zeta=2\int_{\partial D} 
\Phi^a\left(\frac{\zeta+p}{\overline{p}\zeta+1}
\right) \zeta_1\zeta_2 ds_\zeta
\]
Since \(\zeta_1^2+\zeta_2^2=1\), we also get
\[0=\int_{\partial D} \Phi^a\left(\frac{\zeta+p}{\overline{p}\zeta+1}
\right) \zeta_2^2ds_\zeta.
\]
We proceed by induction on \(j\). Suppose that conditions \eqref{cond1ter} 
are valid for \(h=0,...,j\) with \(j<s-1\). We prove that
\begin{eqnarray}\label{cond1jp1}
\int_{\partial D} \Phi^a\left(\frac{\zeta+p}{\overline{p}\zeta+1}
\right) \zeta_1^{j+1-h} \zeta_2^hds_\zeta=0,
\qquad h=0,...,j+1.
\end{eqnarray}
Using the inverse formula (\cite[p.412]{wjc})
\begin{equation}\label{inverseT}
\zeta_1^{j+1}=2^{-j} {\sum^{j+1}_{\substack{{i=0} \\
{j+1-i\>\> even}}}}\!\!\!\!\!\!' \begin{pmatrix}
j+1\\i
\end{pmatrix} T_i(\zeta_1),\quad j\geq 0
\end{equation}
where the prime at the sum symbol means that the first term 
(at \(i=0\)) is to be halved unless it is skipped,
and conditions \eqref{cond1} we get
\[
\int_{\partial D} \Phi^a\left(\frac{\zeta+p}{\overline{p}\zeta+1}
\right) \zeta_1^{j+1}ds_\zeta=0,
\]
that is \eqref{cond1jp1} for \(h=0\). If \(h=2\ell\) with 
\(0<\ell\leq (j+1)/2\), we have, keeping in mind \eqref{inverseT} 
and \eqref{cond1}, 
\begin{eqnarray*}
\int_{\partial D} \Phi^a\left(\frac{\zeta+p}{\overline{p}\zeta+1}
\right) \zeta_1^{j+1-2\ell} \zeta_2^{2\ell}ds_\zeta=
\int_{\partial D} \Phi^a\left(\frac{\zeta+p}{\overline{p}\zeta+1}
\right) \zeta_1^{j+1-2\ell} (1-\zeta_1^2)^\ell ds_\zeta
\\
=\sum_{i=0}^\ell \begin{pmatrix}\ell\\i\end{pmatrix} (-1)^{\ell-i}  
\int_{\partial D}\Phi^a\left(\frac{\zeta+p}{\overline{p}\zeta+1}
\right) \zeta_1^{j+1-2i}  ds_\zeta =0.
\end{eqnarray*}

If \(h=2\ell+1\), with \(0\leq \ell\leq j/2\), we can write 
\eqref{cond1jp1} as
\begin{equation}\label{odd}
\int_{\partial D} \Phi^a\left(\frac{\zeta+p}{\overline{p}\zeta+1}
\right) \zeta_1^{j-2\ell} \zeta_2^{2\ell+1}ds_\zeta
=\sum_{i=0}^\ell \begin{pmatrix}\ell\\i\end{pmatrix} (-1)^{\ell-i}  
\int_{\partial D}\Phi^a\left(\frac{\zeta+p}{\overline{p}\zeta+1}\right)
\zeta_2 \zeta_1^{j-2i} ds_\zeta.
\end{equation}
From the relations \(T_0(\zeta_1)=U_0(\zeta_1)\), 
\(T_1(\zeta_1)=2^{-1}U_1(\zeta_1)\), 
\(T_\kappa(\zeta_1)=2^{-1}(U_\kappa(\zeta_1)-U_{\kappa-2}(\zeta_1))\), 
\(\kappa\geq 2\), and the inverse formula \eqref{inverseT}, we infer 
that \(\zeta_1^{j-2i}\) can be expressed as a linear combination of 
the polynomials \(U_\kappa(\zeta_1)\), \(0\leq\kappa\leq j-2i\). 
Hence, conditions \eqref{cond2} implies that the integrals in the 
right-hand side of \eqref{odd} vanish. \finedim

For example, if \(k=4\) the conditions \eqref{cond1ter} reduce to
\[
\int_{\partial D} \Phi^a\left(\frac{\zeta+p}{\overline{p}\zeta+1}
\right)ds_\zeta=0,\qquad
\int_{\partial D} \Phi^a\left(\frac{\zeta+p}{\overline{p}\zeta+1}\right)
\zeta_jds_\zeta=0, \>\>j=1,2,
\]
which were obtained in  \cite[Proposition 3.12]{lm19}).

As a consequence of Propositions \ref{CS}, \ref{CN} and \ref{CNeS1} 
we deduce the following necessary and sufficient conditions such that 
\(p\in D\) is a point with multipicity \(2s\).
\begin{theorem}\label{CNeS}
Let  \(\Phi=(\phi_1,...,\phi_{2s})\) be an admissible boundary datum and let 
\(\psi_a\) be the solution to \eqref{armonico} 
with boundary datum \(\Phi^a=\sum_{j=1}^{2s}(-1)^j \phi_j\).  
The point  \(p\in{D}\) is a  \(2s-\)point for the function \(U=|\psi_a|\) 
if and only if \(\Phi^a\) satisfies conditions \eqref{cond1ter}.
\end{theorem}

\section{Results on \(6\) species}

In this section we consider  the case of \(6\) competing 
species. As  a consequence of Proposition \ref{pr:formula} 
we show that  \(5\)  limiting configurations are possible. 

\begin{proposition}\label{P3.1}
Let \(U\in\cS\), then only one of the following statement 
is satisfied \\
{\rm i.} $\cZ_3(U)$ consists of one point $a_U\in \overline D$ 
	such that $m(a_U)=6$, \\
{\rm ii.} $\cZ_3(U)$ consists of two points $a_U, b_U\in \overline D$,
	$a_U\neq b_U$, with $m(a_U)=m(b_U)=4$,\\
{\rm iii.} $\cZ_3(U)$ consists of two points 
$a_U, b_U\in \overline D$  with $m(a_U)=3$ and $m(b_U)=5$,\\
{\rm iv.} $\cZ_3(U)$ consists of three different points 
$a_U, b_U,c_U\in \overline D$ with $m(a_U)=4\), 
\(m(b_U)=m(c_U)=3$,\\
{\rm v.} $\cZ_3(U)$ consists of four different points 
$a_U, b_U,c_U,d_U\in \overline D$
with $m(a_U)=m(b_U)=m(c_U)=m(d_U)=3$.
\end{proposition}

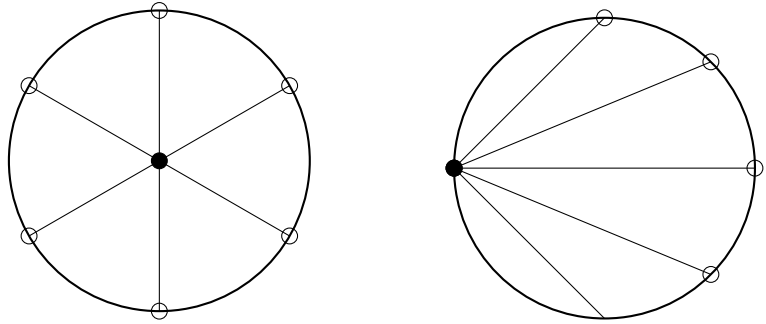
\begin{figure}[h!]
\begin{center}
\begin{tikzpicture}
\draw[thick] (0,0) circle (2cm);
\draw (0,-2) -- (0,2);
\draw (-1.732048,-1) -- (1.732048,1);
\draw (1.732048,-1) -- (-1.732048,1);
\draw (0,-2) circle (3pt);
\draw (0,2) circle (3pt);
\draw (-1.732048,-1) circle (3pt);
\draw (1.732048,1) circle (3pt);
\draw (1.732048,-1) circle (3pt);
\draw (-1.732048,1) circle (3pt);
\filldraw[black] (0,0) circle (3pt);
\end{tikzpicture}
\qquad\qquad
\begin{tikzpicture}
\draw[thick] (0,0) circle (2cm);
\draw (-2,0) -- (0,2);
\draw (-2,0) -- (2,0);
\draw (-2,0) -- (1.414216,1.414216);
\draw (-2,0) -- (1.414216,-1.414216);
\draw (-2,0) -- (0,-2);
\draw (-2,0) circle (3pt);
\draw (-2,0) circle (3pt);
\draw (0,2) circle (3pt);
\draw (2,0) circle (3pt);
\draw (1.414216,1.414216) circle (3pt);
\draw (1.414216,-1.414216) circle (3pt);
\filldraw[black] (-2,0) circle (3pt);
\end{tikzpicture}\end{center}
\caption{Configurations with one 6-point inside \(D\) (on the 
left) and on \(\partial D\) (on the right).}
\label{figtype1Point}
\end{figure}

\begin{figure}[h!]
\begin{center}
\begin{tikzpicture}
\draw[thin] (0,0) circle (2cm);
\draw (-1.414216,1.414216) -- (1.414216,-1.414216);
\draw (-2,0) -- (-0.56,-0.56);
\draw (-2,0) -- (-0.56,-0.56);
\draw (0,-2) -- (-0.56,-0.56);
\draw (2,0) -- (0.56,0.56);
\draw (0,2) -- (0.56,0.56);
\draw (-0.56,-0.56) -- (0.56,0.56);
\draw (-2,0) circle (3pt);
\draw (2,0) circle (3pt);
\draw (0,-2) circle (3pt);
\draw (0,2) circle (3pt);
\draw (-1.414216,1.414216) circle (3pt);
\draw (1.414216, -1.414216) circle (3pt);
\filldraw[black] (-0.56,-0.56) circle (3pt);
\filldraw[black] (0.56,0.56) circle (3pt);
\filldraw[black] (0,0) circle (3pt);
\end{tikzpicture}
\qquad\qquad
\begin{tikzpicture}
\draw[thin] (0,0) circle (2cm);
\draw (0,0) -- (0,0.8);
\draw (0,0) -- (0.69282,-0.4);
\draw (0,0) -- (-0.69282,-0.4);
\draw (0,0.8) -- (1.28172, 1.54);
\draw (0,0.8) -- (-1.28172, 1.54);
\draw (0.69282,-0.4) -- (2.009174,0);
\draw (0.69282,-0.4) -- (0.69282, -1.92);
\draw (-0.69282,-0.4) -- (-2, 0);
\draw (-0.69282,-0.4) -- (-0.69282,-1.92);
\draw (1.28172, 1.54) circle (3pt);
\draw (-1.28172, 1.54) circle (3pt);
\draw (-2,0) circle (3pt);
\draw (-0.69282,-1.92) circle (3pt);
\draw (2.009176,0) circle (3pt);
\draw (0.69282, -1.92) circle (3pt);
\filldraw[black] (0.69282,-0.4) circle (3pt);
\filldraw[black] (-0.69282,-0.4) circle (3pt);
\filldraw[black] (0,0.8) circle (3pt);
\filldraw[black] (0,0) circle (3pt);
\end{tikzpicture}

\end{center}
\caption{Configurations when $\cZ_3(U)$ consists of three 
points  (on the left) and four points (on the right).}
\label{figtype3Points}
\end{figure}
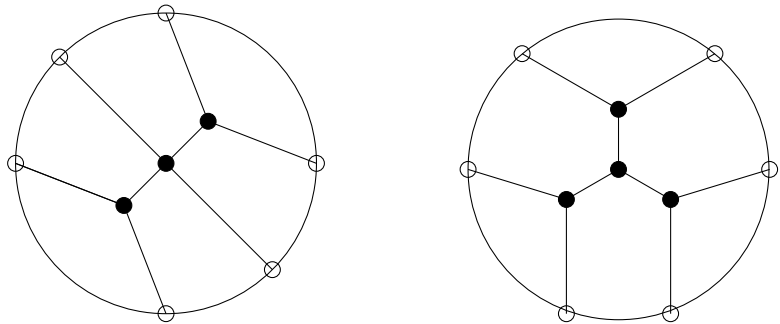

\dimo The set \(\cZ_3(U)\) is nonempty and contains at 
most \(4\) points. If \(\cZ_3(U)=\{a_U\}\) then, from 
\eqref{formula}, we deduce that \(m(a_U)=6\) (cf. figure 
\ref{figtype1Point}). If \(\cZ_3(U)=\{a_U,b_U\}\) then, 
from \eqref{formula}, we deduce that \(m(a_U)+m(b_U)=8\), 
that is ii. or iii. (cf. figure \ref{figtype2Points}).
If \(\cZ_3(U)=\{a_U,b_U,c_U\}\) then, from \eqref{formula}, 
we deduce that \(m(a_U)+m(b_U)+m(c_U)=10\), that is iv. 
(cf. figure \ref{figtype3Points} on the left).
If \(\cZ_3(U)=\{a_U,b_U,c_U,d_U\}\) then, from \eqref{formula}, 
we deduce that \(m(a_U)+m(b_U)+m(c_U)+m(d_U)=12\), that is v. 
(cf. figure \ref{figtype3Points} on the right).

\begin{figure}[h!]
\begin{center}
\begin{tikzpicture}
\draw[thick] (0,0) circle (2cm);
\draw (0,2) -- (0,-2);
\draw (-1.618032,1.175568) -- (1.618032,1.175568);
\draw (-1.618032,-1.175568) -- (1.618032,-1.175568);
\draw (1.618032,1.175568) circle (3pt);
\draw (1.618032,-1.175568) circle (3pt);
\draw (-1.618032,1.175568) circle (3pt);
\draw (-1.618032,-1.175568) circle (3pt);
\draw (0,2) circle (3pt);
\draw (0,-2) circle (3pt);
\filldraw[black] (0,1.175568) circle (3pt);
\filldraw[black] (0,-1.175568) circle (3pt);
\end{tikzpicture}
\qquad\qquad
\begin{tikzpicture}
\draw[thick] (0,0) circle (2cm);
\draw (0,0) -- (1.902112,-0.6180336);
\draw (0,0) -- (1.175568,1.618032);
\draw (0,0) -- (-1.175568,1.618032);
\draw (0,0) -- (-1.902112,-0.6180336);
\draw (0,0) -- (0,-0.8);
\draw (0,-0.8) -- (1.302504, -1.552);
\draw (0,-0.8) -- (-1.302504, -1.552);
\draw (1.902112,-0.6180336) circle (3pt);
\draw (1.175568,1.618032) circle (3pt);
\draw (-1.175568,1.618032) circle (3pt);
\draw (-1.902112,-0.6180336) circle (3pt);
\draw (1.302504,-1.552) circle (3pt);
\draw (-1.302504,-1.552) circle (3pt);
\filldraw[black] (0,0) circle (3pt);
\filldraw[black] (0,-0.8) circle (3pt);
\end{tikzpicture}
\end{center}
\caption{Configurations when $\cZ_3(U)$ consists of two 
points $a_U\neq b_U$, with $m(a_U)=m(b_U)=4$ (on the left) 
and $m(a_U)=3$ and $m(b_U)=5$ (on the right).}
\label{figtype2Points}
\end{figure}
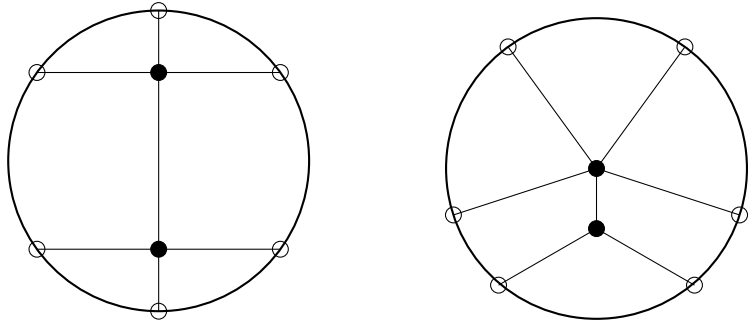 

\finedim
 
In the case of \(6\) species, Theorem \ref{CNeS} can be formulated as follows. 

\begin{proposition}\label{6pointCS}
Let \(\Phi=(\phi_1,...,\phi_{6})\)  be an admissible datum, 
let \(U\in\cS\) and let \(\psi_a\) be the solution to 
\eqref{armonico} with boundary datum 
\(\Phi^a=\sum_{j=1}^{6}(-1)^j \phi_j\).  The function \(U\) 
has a 6-point \(p\in{D}\) if and only if the boundary datum 
\(\Phi=(\phi_1,...,\phi_6)\)  satisfies the conditions
\begin{eqnarray}\label{C1}
&\displaystyle\sum_{j=1}^6 (-1)^j \int_{\partial D} \phi_j\left(
	\frac{\zeta+p}{\overline{p}\, \zeta+1}\right) ds_\zeta&=0\,,\\\label{C2}
&\displaystyle\sum_{j=1}^6 (-1)^j \int_{\partial D} \phi_j\left(
	\frac{\zeta+p}{\overline{p}\, \zeta+1}\right) \zeta_r ds_\zeta
	&=0\,, \qquad r=1,2
\end{eqnarray}
	\begin{eqnarray}\label{C3}
	&\displaystyle\sum_{j=1}^6 (-1)^j \int_{\partial D} \phi_j\left(
	\frac{\zeta+p}{\overline{p}\, \zeta+1}\right) \zeta_1^2ds_\zeta&=0\,,\\\label{C4}
&\displaystyle\sum_{j=1}^6 (-1)^j \int_{\partial D} \phi_j\left(
	\frac{\zeta+p}{\overline{p}\, \zeta+1}\right) \zeta_1\zeta_2 ds_\zeta
	&=0.
\end{eqnarray}
 Moreover the function \(U=|\psi_a|\).
\end{proposition}
\dimo
Conditions \eqref{C1}-\eqref{C4} are easily obtained by assuming 
\(s=3\) in \eqref{cond1ter}. \finedim

\begin{proposition}\label{second}
Let \(\Phi\) be an admissible datum and let \(\psi_a\) be the 
solution to \eqref{armonico} 
with boundary datum \(\Phi^a=\sum_{j=1}^{6}(-1)^j \phi_j\). 
Conditions \eqref{C1}-\eqref{C4} are equivalent to
\begin{equation}\label{101}
\psi_a(p)=0,\quad \nabla \psi_a(p)=(0,0),\quad H\psi_a(p)=\bf 0
\end{equation}
where \(H\psi_a=\{\partial^2_{x_ix_j} \psi_a\}_{i,j=1,2}\) 
denotes the hessian matrix of the function \(\psi_a\).
\end{proposition}

\dimo
We introduce the transformation \eqref{conforme} and define 
\(\widetilde \Psi_a\) and \(\widetilde \Phi^a\) according to \eqref{TR}.
We have
\[
\psi_a(p)=\psi_a(R_p(0))=\widetilde \Psi_a(0);
\]
\[
\nabla_x \psi_a(p)=\begin{pmatrix}
(1-|p|^2)^{-1}&0\\
0&(1-|p|^2)^{-1}
\end{pmatrix} \nabla_\zeta\widetilde \Psi_a(0)
\]
and, assuming \(p=(p_1,p_2)\), 
\[
H_x\psi_a(p)=\frac{1}{(1-|p|^2)^2} H_\zeta \widetilde\Psi_a(0)
+\frac{2}{1-|p|^2} \left(\partial_{\zeta_1} \widetilde\Psi_a(0) 
\begin{pmatrix} p_1&p_2\\p_2&-p_1
\end{pmatrix}+\partial_{\zeta_2} \widetilde\Psi_a(0) 
\begin{pmatrix} -p_2&p_1\\p_1&p_2
\end{pmatrix}\right)\,.
\]
By construction, \(\widetilde \Psi_a\) solves \eqref{DP2}. 
By the Poisson integral formula \eqref{poisson} we deduce that
\[
\widetilde \Psi_a(0)=\frac{1}{2\pi}\int_{\partial D} 
\frac{\widetilde \Phi^a (\eta)}{|\eta|^2}ds_\eta=
\frac{1}{2\pi}\int_{\partial D} \widetilde \Phi^a (\eta)ds_\eta
= \frac{1}{2\pi}\int_{\partial D}  \Phi^a (R_p(\eta))ds_\eta.
\]
On the other hand, by direct differentation of the Poisson integral,  
we get
\[
\frac{\partial}{\partial \zeta_j} \widetilde\Psi_a (0)
= \frac{1}{\pi} \int_{\partial D} \widetilde \Phi^a(\eta)\eta_j ds_\eta 
= \frac{1}{\pi} \int_{\partial D}  \Phi^a(R_p(\eta))\eta_j ds_\eta,
\quad j=1,2
\]
\[
\frac{\partial^2}{\partial^2 \zeta_j} \widetilde\Psi_a (0)
= -\frac{2}{\pi} \int_{\partial D} \widetilde \Phi^a(\eta)ds_\eta+
\frac{4}{\pi} \int_{\partial D} \widetilde \Phi^a(\eta)\eta_j^2 ds_\eta,
\qquad j=1,2
\]
\[
\frac{\partial^2}{\partial \zeta_1\partial \zeta_2} \widetilde\Psi_a (0)
=\frac{4}{\pi} \int_{\partial D} \widetilde \Phi^a(\eta)\eta_1 \eta_2 ds_\eta.
\]
The equivalence between \eqref{C1}-\eqref{C4} and \eqref{101} easily follows.
\finedim

\begin{proposition}\label{6pointBIS}
Let \(\Phi=(\phi_1,...,\phi_{6})\) be an admissible datum and suppose 
that the related harmonic function \(\psi_a\), solution to \eqref{armonico} 
with boundary datum \(\phi^a\), has two critical points \(p,q\in D\) such 
that \(\psi_a(p)=\psi_a(q)=0\), then \(p,q\) are 4-points for the 
function \(U=|\psi_a|\in\cS\).
\end{proposition}

\dimo By standard theory of harmonic functions the zero set of 
\(\psi_a\) around a critical point at level 0 is made by (at least) 
\(4\) half-lines, meeting with equal angles. We infer that locally 
around \(p\)  the function \(\psi_a\) defines \(k_p\) nodal components 
with \(k_p\geq 4\) and \(k_p\) is even because \(\psi_a\) has alternate 
positive or negative sign on adjacent sets. 
Analogously around  \(q\) the function \(\psi_a\) defines \(k_q\) 
nodal components with \(k_q\geq 4\) and \(k_q\) is even. 
By Proposition \ref{P3.1} we infer that \(k_p=k_q=4\). 
Hence the function \(U=|\psi_a|\) is nonnegative, satisfies the boundary 
datum \(\Phi\), has exactly \(6\) nodal regions and
generates an element of \(\cS\), with \(p,q\in \cZ_3(U)\).
\finedim

Conversely, as a direct application of Proposition \ref{iffharmonic2} 
we have the following result.
\begin{proposition}\label{6pointTER}
Let \(\Phi=(\phi_1,...,\phi_{6})\) be an admissible datum. 
Let \(U\in \cS\) such that \(\cZ_3(U)\) contains 
two points \(p,q\) in \({D}\) with  multiplicity \(4\). 
Then  \(U=|\psi_a|\), where \(\psi_a\) is the harmonic 
function such that \(\psi_a=\Phi^a= \sum_{j=1}^{6} (-1)^j\phi_j\) 
on \(\partial D\). 
Moreover \(p,q\) are critical points for \(\psi_a\) at zero level.
\end{proposition}

\section{Conclusion and open problems}

Propositions \ref{6pointBIS} and \ref{6pointTER} show that 
the limiting configurations, which elements of \(\cZ_3(U)\) 
have only even multiplicity, are closely connected to harmonic 
solutions of \eqref{armonico} with alternate boundary datum.
Since such solutions have to satisfy some integral conditions
(see \eqref{cond1ter}), it follows that the most probable
segregation configurations have only points in \(\cZ_3(U)\)
with odd multiplicity (see also Remark 3.13 in \cite{lm19}).

In fact we think that the most probable configurations in nature 
are those with points with multiplicity three.
But this means that odd multiple points (with multiplicity
greater than 3) have to verify some conditions  
in order to belong to subsets of admissible data with nontrivial 
codimension. Now this is an open problem, since our techniques
are strongly related to the properties of the harmonic functions,
and these functions are related only to 
critical points with even multiplicity.

An open problem we wish to handle is to compare results by 
reaction-diffusion systems in populations competition, as in 
\cite{bz,bz2,ctv2,ctv4,dwz,lm19,tvz}, and results obtained 
by differential games theory contained in \cite{mk05,halw} 
(see also the references therein), in order to reach a better 
understanding of the territoriality of the competing 
species or groups (see also \cite{bdf08}).

\subsection*{Acknowledgment}
The authors are very grateful to the anonymous referees 
for their insightful comments and helpful suggestions.

\end{document}